\documentstyle[12pt]{article}

\scrollmode

\renewcommand{\a}{\alpha}         
\newcommand{\De}{\Delta}          
\newcommand{\de}{\delta}          
\newcommand{\Ga}{\Gamma}          
\newcommand{\ga}{\gamma}          
\newcommand{\Om}{\Omega}          
\newcommand{\om}{\omega}          
\newcommand{\Th}{\Theta}          
\renewcommand{\th}{\theta}        

\newcommand{\A}{{\cal A}}         
\newcommand{\ENE}{{\cal N}}       
\newcommand{\M}{{\cal M}}         
\newcommand{\Tau}{{\cal T}}       
\newcommand{\X}{{\cal X}}         

\renewcommand{\S}{{\bf S}}        
\newcommand{\T}{{\bf T}}          

\newcommand{\opname}[1]{\mathop{\rm#1}\nolimits} 
\newcommand{\id}{\opname{id}} 

\newcommand{\Tk}{{{\bf T}^{(k)}}}     
\newcommand{\Tkmu}{{{\bf T}^{(k-1)}}} 
\newcommand{\Tu}{{{\bf T}^{(1)}}}     
\renewcommand{\v}{\vee}               
\newcommand{\Xk}{{X^{(k)}}}           
\newcommand{\Xkmu}{{X^{(k-1)}}}       
\newcommand{\Xl}{{X^{(l)}}}           
\renewcommand{\:}{\colon}         

\def\<#1>{\langle#1\rangle}       

\newcommand{\pd}[2]{\frac{\partial#1}{\partial#2}} 

\newcommand{\sepword}[1]{\qquad\mbox{#1}\quad} 
 
\newtheorem{thm}{Theorem}

\begin{document}

\begin{center}{\bf {Higher order Lagrangian supermechanics}}

\bigskip

{\it {Jos\'e F. Cari\~nena$^\dag$\  and\ 
 H\'ector Figueroa$^\ddag$}}

 \medskip

${}^{\dag}$ Departamento de F\'{\i}sica Te\'orica,
            Universidad de Zaragoza,
            
 50009 Zaragoza, Spain

 \medskip
 
${}^\ddag$ Escuela de Matem\'atica, 
            Universidad de Costa Rica,
            
 2060 San Jos\'e, Costa Rica
 \end{center}



Generalized mechanics where the Lagrangian 
describing the system depends on higher order 
derivatives has been investigated in connection
with relativistic dynamics of particles, supersymmetry, 
polymer physics, string models and gravity theory, among 
other problems (see, e.g., \cite{Zi}).
 For instance, it has been proposed to add to the point 
particle action a term proportional to acceleration,
or the inclusion of higher order terms in
superstring action. The mathematical theory of such 
generalized Lagrangians was started by Ostrogradski
\cite{Ostrogradski}, and it was more recently
analysed, from a geometric perspective, in \cite{Pepinv}. 
The study of classical higher order Lagrangian 
systems was considered in 
\cite{CrampinSarlet}. Our aim is to provide the
appropriate mathematical framework to deal with such systems,
along similar lines to those developed in~\cite{Hector}.

  The first order Lagrangian formalism of supermechanics
takes place on the tangent supermanifold $T\M$ of a given 
graded manifold $\M=(M,\A)$, \cite{Hector,Ibort}. 
As in the non--graded case,  using the geometric structure
of the tangent supermanifold, which is encoded in the
vertical superendomorphism $S$, one can construct a graded
symplectic structure $\Om_L=-d\Th_L=-d(dL\circ S)$, 
out of a regular 
Lagrangian superfunction $L\in T\M$. The dynamics of
such systems is governed by a second order differential
superequation $\Ga\in\X(T\M)$, which is the unique 
solution of the dynamical equation $i_X\Om_L=dE_L$, where
the superenergy is defined, in terms of the Liouville
supervector field $\De$, by $E_L:=\De L-L$.

  In analogy to the classical case, to establish a 
one--to--one correspondence between symmetries of $L$
and constants of motion it is necessary to generalize
the concept of symmetry so as to include ``non--point
transformations". These symmetries turn out NOT to be 
supervector fields but rather supervector fields along
a morphism (see \cite{Hectoriv} for a detailed
discussion of supervector fields and graded forms
along a morphism). 
  Moreover, from a technical point of view, it turns out 
to be convenient to extend such supervector fields so 
one can work in the second order tangent supermanifold
$T^2\M$ \cite{Hector}. In particular, it is useful 
to rewrite 
the dynamical equation as $\ga^*(\de L)=0$, where
$\ga$ is a section associated to $\Ga$ \cite{Hector}
(we shall describe $\ga$ in the general case below),
and $\de L$ is the Euler--Lagrange graded 1--form
\begin{equation}
\de L:= i_{\Tu}\Om_L - \Tau^*_{2,1}d\,E_L.
\label{uno}
\end{equation}
Here $\Tau_{2,1}\:T^2\M\to T\M$ denotes the canonical 
projection, and $\Tu$ is the time derivative operator,
which is a canonical supervector field
along $\Tau_{2,1}$ to be described later.

   We see that already in first order Lagrangian
supermechanics higher order tangent supermanifolds have
a role to play. The most intuitive way to define the
$k$--th order tangent manifold of $M$, in classical 
geometry, is as the space of equivalence classes of 
curves on $M$ that agree up to order $k$ 
\cite{CrampinSarlet,Pepinii}.  Nevertheless,
this approach is not easy to generalize to graded
manifolds since it is based on a construction using 
points,
and the information of a graded manifold is not 
concentrated on the underlying manifold but rather on 
the sheaf of superalgebras. Thus, in the graded context,
it is better to define the $k$--th order tangent 
supermanifold recursively using a ``diagonal process'',
so for instance to define $T^2\M$ one takes advantage of
the two natural projections of $T(T\M)$ onto $T\M$,
namely the tangent projection of $T(T\M)$, and the 
tangent morphism of the projection $T\M\to \M$. To
define $T^3\M$ one replaces $T\M$ by $T^2\M$ in the
above construction and so on (see \cite{Hectorv}
for details).

   The most important technical tool when dealing 
with higher order tangent supermanifolds is the
total time derivative operator.  The point is that
any graded 
manifold $\ENE$ has associated a canonical supervector
field along its tangent projection $T\ENE\to\ENE$ 
\cite{Hectoriv}. When 
the graded manifold is $\ENE=T^k\M=(T^kM,T^k\A)$ we 
obtain a supervector field $\tilde\T^k\:T^k\A\to T(T^k\A)$
along the tangent projection of $T^k\M$. The {\it total
time derivative of order} $k+1$ is the supervector
field $\Tk\:T^k\A\to T^{k+1}\A$ along the canonical 
projection $\Tau_{k+1,k}\:T^{k+1}\M\to T^k\M$ defined by
$\Tk:= i^*_{k+1}\circ\tilde\T^k$, where 
$i_{k+1}\:T^{k+1}\M\to T(T^k\M)$ is the canonical inclusion.

   Total time derivatives are the brigdes that
allow us to lift objects from one tangent supermanifold
to another one of higher order. For instance, out of a
superfunction $f$ in $\A$ we can construct $k+1$ 
superfunctions in $T^k\A$, to wit
$f^k_j := \tau_{k,j}^*\circ\T^{(j-1)}\circ\dots
\circ\Tu\circ\T(f)$. In particular, when $f$ runs through
a local system of supercoordinates on $\M$, we obtain
local supercoordinates for $T^k\A$. Thus, 
we can also define the {\it lift of a supervector field}
$X$ along the projection $\Tau_{k,0}$ as the
supervector field $\Xl$ along the projection 
$\Tau_{k+l,l}$ given by the equations:
$\Xl(f^l_j):= i^*_{k,l}\bigl((Xf)^l_j\bigr)$ for all 
$f\in\A$, where again $i_{k,l}$ denotes the canonical
inclusion.

  As is well known, the geometrical
information of a cotangent bundle is concentrated in its
canonical symplectic form. In the same way, the 
geometrical information of
a tangent bundle (of any order) is centered in its 
vertical endomorphism and its Liouville vector field.  
To define these objects in
the context we are interested in, we start by lifting
superfunctions adapting some ideas of Tulczyjew 
\cite{Tulczyjewi}, (see also
\cite{Pepinii,CrampinSarlet}). Thus, the {\it vertical
lift of a superfunction} $f\in T^{k-1}\A$ is the super
function $f^V\in T^k\A$ defined by
$$
f^V := 
\sum^m_{i=1}\sum^{k-1}_{j=0} 
{1\over j+1}\pd{F}{q^i_{j,k}}q^i_{j+1,k} 
+ \sum^n_{\a=1}\sum^{k-1}_{j=0}
{1\over j+1} \pd{F}{\th^\a_{j,k}}\th^\a_{j+1,k}, 
$$
where $F:=\tau^*_{k,k-1}f$, and $(q^i_{j,k},\th^\a_{j,k})$
are the supercoordinates on $T^k\M$ obtained from a
local system of supercoordinates $(q^i,\th^\a)$
on $\M$. It turns out that, the
action of a supervector field on these superfunctions
determines completely the  supervector field, so we
can define the {\it vertical lift of a supervector field
along} $\Tau_{k,k-1}$ as the supervector field on $T^k\M$
that satisfies the equations $X^V(f^V) = X(f)$, 
for all $f$ in $T^{k-1}\A$. The {\it Liouville 
supervector field} on $T^k\M$ is nothing but the vertical
lift of the total time derivative operator, in other 
words, is the supervector field $\De_k$ on $T^k\M$ 
defined by $\De_k:=(\Tkmu)^V$.
In this sense, we noticed that the time derivative
operator is a more fundamental object than the Liouville
supervector field.
 
On the other hand, the graded
tensor field on $T^k\M$ of type $(1,1)$ 
$S_k \: \X(T^k\A) \to \X(T^k\A)$,  defined by
$S_k(Y) :=(Y\circ\tau^*_{k,k-1})^V$,
is called the {\it vertical superendomorphism}.
As in the case $k=1$, given a regular super Lagrangian 
$L$, which now is a superfunction in $T^k\M$, 
the symplectic structure of the dynamical system
associated to $L$ is
built out of the vertical superendomorphism, but in a
somewhat more intricate way.  With this in mind, we
consider the {\it Cartan operator}
$\S^{(k)}\:\Om^1(T^k\A)\to \Om^1(T^{2k-1}\A)$ 
defined by the formula
$$
\S^{(k)}:=\sum^k_{l=1}{(-1)^{l+1}\over l!}
\Tau^*_{2k-1,k+l-1}\circ d^{l-1}_{\T^{(k)}}\circ 
S^{*l}_k,
$$
where $d^l_{\T^{(k)}}:=d_{\T^{(k+l-1)}}\circ 
d_{\T^{(k+l-2)}}\circ\cdots\circ d_{\T^{(k)}}$ and
$d_{\T^{(r)}}$ denotes the $\Tau^*_{r+1,r}$--derivation 
that extends $\T^{(r)}$ to the full Cartan algebra of 
graded forms \cite{Hector}, which is the analogue of
a Lie derivative for this kind of supervector fields.

   We now introduce the usual graded forms that are
used in the description of Lagrangian formalisms.  The
{\it Cartan 1--form} and the {\it Cartan 2--form}
associated to a super Lagrangian $L$ are, respectively,
the graded forms on $T^{2k-1}\M$ given by
$$
\Th_L:=\S^{(k)}(d\,L) \sepword{and}
\Om_L:= -d\,\Th_L.
$$
On the other hand, the {\it Euler--Lagrange} form 
associated to $L$ is defined by 
\begin{equation}
\de L= \Tau^*_{2k,k}(d\,L)- d_{\T^{(2k-1)}}\Th_L.
\label{dos}
\end{equation}
We notice that $\Th_L$ and
$\de L$ are $\Tau_{2k-1,k-1}$--semibasic and 
$\Tau_{2k,0}$--semibasic, respectively. The importance
of being $\Phi$--semibasic, where $\Phi\:\ENE\to\M$ is
a morphism, 
lies in the fact that there
is a one--to--one correspondence between these graded
1--forms on $\ENE$ and graded 1--forms along $\Phi$ 
\cite{Hector}: if $\om$ is a $\Phi$--semibasic graded 
1--form on $\ENE$, the corresponding  graded 1--form 
along $\Phi$ is the one defined by $\check\om(X):=\om(Y)$,
whenever $X\in\X(\Phi)$ is such that $X=Y\circ\phi^*$ for
some supervector field $Y$ on $\ENE$. Moreover, one
can check that if $\om$ is a $\Tau_{k,l}$--semibasic 
graded 1--form on $T^k\M$, and $X$ is a
supervector field along $\Tau_{r,0}$, then
\begin{equation}
i_{\Xk}\om = \tau^*_{k+r,k}\bigl(
\<\tau^*_{k,l+r}\circ\Xl,\check\om>\bigr),
\label{tres}
\end{equation}
provided that $r+l\leq k$ and $l\leq r$. Thus, using the
Cartan identity and (\ref{tres}) 
\begin{eqnarray*}
\de L & = & \Tau^*_{2k,k}(dL)-i_{\T^{(2k-1)}}d\,\Th_L 
              - d\,i_{\T^{(2k-1)}}\Th_L  \\
 & = & i_{\T^{(2k-1)}}\Om_L 
              - d\Bigl(i_{\T^{(2k-1)}}\Th_L-\tau^*_{2k,k}L\Bigr) \\
 & = & i_{\T^{(2k-1)}}\Om_L - \Tau^*_{2k,2k-1}
        d\Bigl( \<\tau^*_{2k-1,k}\circ\Tkmu,\Th^\v_L>
                - \tau^*_{2k-1,k}L \Bigr).
\end{eqnarray*}
A comparison of this equation and (\ref{uno}) leads us to
define the {\it energy} associated to
the super--Lagrangian $L$ as the superfunction in
$T^{2k-1}\A$ given by
$$
E_L:= \<\tau^*_{2k-1,k}\circ\Tkmu,\Th^\v_L>
-\tau^*_{2k-1,k}L.
$$

  In this way, we associate to a super Lagrangian 
function $L$ the Hamiltonian dynamical system 
$(T^k\M,\Om_L,E_L)$. The dynamics of such 
systems is governed by a superdifferential equation
of order $k+1$, which, by definition, is a supervector
field on $T^k\M$ that satisfy either of the equivalent
equations
$$
\Ga \circ \tau^*_{k,k-1} = \Tkmu
\sepword{or}
S_k(\Ga)=\De_k.
$$
First thing to notice is that super differential
equations are always even supervector fields. On
the other hand, they can also be regarded, at least 
formaly, as sections of a bundle,
since there is a one--to--one correspondence
between super differential equations of order $k+1$
and sections of $\Tau_{k+1,k}\:T^{k+1}\M\to T^k\M$, 
in the sense that each
super differential equations $\Ga$ of order $k+1$ is
associated to a morphism $\ga\: T^{k+1}\A \to T^k\A$,
such that $\ga \circ\tau^*_{k+1,k} = \id_{T^k\A}$, and
viceversa.

  Thus, it is only natural to say that a supervector
field $\Ga$ is {\it Lagrangian with respect} to a
super Lagrangian $L\in T^k\M$ if it is a superdifferential
equation of order $2k$ and satisfies the {\it dynamical
equation}
$$
i_\Ga\Om_L = d\,E_L.
$$
As in the classical case, when $L$ is regular, that is
when $\Om_L$ is non--degenerated, there exists a unique
supervector field satisfying the dynamical equation.
Moreover, this equation can be reformulated as 
$\ga\circ \de L=0$, where $\ga$ is the section associated
to $\Ga$.

  One can establish a relation between constants of motion 
and supervector fields:

\begin{thm}
	Let $L\in T^k\A$ be a regular super--Lagrangian. A 
	superfunction $G \in T^{2k-1}\A$ is a constant of motion 
	if, and only if, there exists $X \in \X(\Tau_{2k-1,0})$ 
	such that 
	$$
	\T^{(2k-1)} G = -\<\tau^*_{2k,2k-1} \circ X,(\de L)^\v>.
	$$
\end{thm}

  Because of this theorem, and taking in consideration
the original definition in the classical context
given by Marmo and Mukunda \cite{Marmo} (see also
\cite{Pepinii}), we say that a supervector field $X$
along the projection $\Tau_{2k-1,0}$ is a
{\it generalized infinitesimal supersymmetry} of the 
dynamical system $(T^k\M,\Om_L,E_L)$, if there exists a 
superfunction $F\in T^{3k-2}\M$ such that
\begin{equation}
\Xk L = \T^{(3k-2)} F.
\label{cuatro}
\end{equation}
Moreover, we can prove a generalization
of Noether theorem and its converse \cite{Hectorv}:

\begin{thm}
	Let $X$ be a generalized infinitesimal supersymmetry of a
	regular super Lagrangian $L$. Then there exists 
	a superfunction $G\in T^{2k-1}\A$ such that is a constant 
	of motion and satisfies $\tau^*_{3k-1,2k-1}G= \<\Xkmu,
	(\Tau^*_{3k-2,2k-1}\Th_L)^\v> -F$,
	where $F$ is the superfunction associated to $X$ via
	(\ref{cuatro}). 
	
	Reciprocally, if $G$ is a constant of motion of the
	dynamical system associated to $L$, then there exists a 
	supervector field $X$ along $\Tau_{2k-1,0}$ such that 
	the superfunction 
	$F:=\<\Xkmu,(\Tau^*_{3k-2,2k-1}\Th_L)^\v> - 
	\tau^*_{3k-1,2k-1}G$ satisfies (\ref{cuatro}).
\end{thm}

  So we obtain an explicit one--to--one correspondence 
between constants of motion and generalized infinitesimal
supersymmetries of the dynamical system 
$(T^k\M,\Om_L,E_L)$.

\end{document}